\documentclass{article}
\usepackage{geometry}
\usepackage[utf8]{inputenc}
\usepackage{amsmath}
\usepackage{amssymb}
\usepackage{amsthm}
\usepackage{graphicx}
\usepackage{tikz-cd}
\usepackage{url}

\theoremstyle{definition}
\newtheorem{example}{Example}[section]
\newtheorem{remark}{Remark}[section]

\theoremstyle{plain}
\newtheorem{lemma}{Lemma}[section]
\newtheorem{theorem}{Theorem}[section]
\newtheorem{corollary}{Corollary}[section]

\numberwithin{equation}{section}

\newcommand{\C}{\mathbb{C}}
\renewcommand{\P}{\mathbb{P}}
\newcommand{\E}{\mathbb{E}}
\newcommand{\T}{\mathcal{T}}
\newcommand{\Var}{\text{Var}}

\newcommand{\lt}{|L(T)|}
\newcommand{\p}{P(T)}
\newcommand{\pt}{P_E(T)}
\newcommand{\w}{W(T)}
\newcommand{\wt}{W_E(T)}
\newcommand{\q}{Q(T)}
\newcommand{\qt}{Q_E(T)}
\newcommand{\bt}{\binom{|T|}{2}}
\newcommand{\blt}{\binom{\lt}{2}}

\newcommand{\Gr}{G_{\text{r}}}
\newcommand{\Gd}{G_{\text{$d$}}}
\newcommand{\Gp}{G_{\text{$\alpha$}}}

\newcommand{\Fr}{F_{\text{r}}}
\newcommand{\Fd}{F_{\text{$d$}}}
\newcommand{\Fp}{F_{\text{$\alpha$}}}

\title{The External Path Length and Terminal Wiener Index in Random Trees\footnote{An extended abstract of this paper was presented at the 37th International Conference on Probabilistic, Combinatorial and Asymptotic Methods for the Analysis of Algorithms (AofA26), see \cite{ext_abs}.}}

\usepackage{authblk}
\usepackage{orcidlink}

\author[1]{Jacob Lundblad \orcidlink{0009-0003-0032-5905}\thanks{Supported by the Wallenberg AI, Autonomous Systems and Software Program (WASP) funded by the Knut and Alice Wallenberg Foundation.}}
\author[1,2]{Stephan Wagner \orcidlink{0000-0001-5533-2764}\thanks{Supported by the Swedish research council (VR), grant 2022-04030.}}
\affil[1]{Department of Mathematics, Uppsala University, Sweden} 
\affil[2]{Institute of Discrete Mathematics, TU Graz, Austria}

\date{}

\begin{document}

\maketitle

\begin{abstract}
We consider two distance-based concepts in rooted trees and compare them to their leaf-restricted counterparts. The path length is the sum of all distances from the root to the vertices of the tree, while the external path length is the sum of all distances from the root to the leaves of the tree. Upon dividing by the number of vertices and leaves respectively, we obtain the average distance to the root. The Wiener index is the sum of all distances between pairs of vertices, and the terminal Wiener index is the sum of all distances between pairs of leaves. Dividing by the number of pairs of vertices and pairs of leaves respectively, we obtain the average distance between vertices and average distance between leaves. For two important classes of random trees, we show that the average distance of the root to a random leaf is almost the same as the average distance to a random vertex and that the average distance between vertices is almost the same as the average distance between leaves. For Bienaym\'e--Galton--Watson trees, each of these differences is bounded in probability. For three varieties of increasing trees (recursive trees, $d$-ary increasing trees, generalised plane-oriented recursive trees), on the other hand, the differences converge to constants in probability.
This allows the transfer of convergence results for the path length to the external path length, and results for the Wiener index to be transferred to the terminal Wiener index.

\end{abstract}

\section{Introduction}

The \emph{path length}, defined as the sum of the depths (distances to the root), is a fundamental quantity associated with a rooted tree. Different notions of path length exist depending on the range of summation: the \emph{internal path length} is the sum of the depths of all internal vertices (non-leaves), while the \emph{external path length} is the sum of the depths of all external vertices (leaves). When the sum is taken over all vertices, one obtains the \emph{(total) path length} (also known as \emph{total height}) of the tree.

The Wiener index of a graph, defined as the sum of pairwise distances between vertices, was first introduced by chemist Harold Wiener in 1947 in relation to chemical properties of molecules (in particular the boiling points of alkane molecules)~\cite{HaroldWiener}. Restricting this sum to pairs of leaves one instead obtains the \emph{terminal Wiener index} (external Wiener index), which was first introduced by Gutman, Furtula and Petrovi\'c \cite{TerminalWienerIndexGutman}.

The path length and Wiener index in their different versions can be used as a measure for the performance of algorithms based on trees. As a result, their asymptotic distributions have been studied for various types of random trees.

In this paper, we will be particularly interested in how the different notions of path length and Wiener index differ. If we divide the internal/external/total path length of a tree by the internal/external/total number of vertices, we obtain the respective average distance to the root. It is intuitive to assume that the average distance from a leaf to the root should be greater than the average distance from an arbitrary vertex to the root (though this is not deterministically the case: consider for example a tree with $n$ leaves and a path of length $n$ attached to a common root as a counterexample). However, the question remains how much these averages differ.
Analogously, dividing the Wiener index and terminal Wiener index by the number of pairs of vertices/leaves yields the average distances between the vertices/leaves. Again it is intuitive to assume that the leaves are on average further from each other than general vertices, but this is again not deterministically true, and how much these averages differ remains to be determined.

We study these questions for two different types of random trees---conditioned Bienaym\'e--Galton--Watson trees (critical, aperiodic, with finite variance) and linear preferential attachment trees (which encompass recursive trees, $d$-ary increasing trees and generalised plane-oriented recursive trees). In both cases, we will show that the difference between the average distance between two arbitrary vertices and the average distance between two leaves, as well as the difference between the average depth of a leaf and the average depth of an arbitrary vertex, are both only of constant order. This implies in particular that the different notions of Wiener index lead to the same limiting distribution, something that also holds for the different notions of path length.

In particular, for both types of trees, we obtain the asymptotic behaviour of the first and second moments of the difference between the average distance from the root to a random vertex and the average distance from the root to a random leaf, as well as for the difference between the average distance between vertices and the average distance between leaves. In all four cases, these are asymptotically constant. For linear preferential attachment trees, the second moments are asymptotically equal to the squares of the first, so that the variances go to zero and we obtain convergence in probability to constants. For Bienaym\'e--Galton--Watson trees, we do not get quite as strong results, but the differences are shown to be bounded in probability.

For all these families, the path length and related depth parameters have previously been studied individually.
For conditioned Bienaym\'e--Galton--Watson trees, it is well known that the expectation of the path length is of asymptotic order $n^{3/2}$. After a proper rescaling, it converges weakly to the \emph{Airy distribution}, which can be seen as the area under a Brownian excursion. This was first shown by Tak\'acs (in special cases, \cite{Takacs1992,Takacs1993}) and Aldous (for more general Bienaym\'e--Galton--Watson trees, see \cite[Corollary 9]{Aldous1991}) in the 1990s, see also \cite[Theorem 4.9]{Drmota2009}.

On the other hand, for trees of logarithmic height (such as linear preferential attachment trees), the path length is typically of order $n \log n$. The expected value in the case of recursive trees has the nice closed form $n(H_n-1)$, where $H_n$ denotes the $n$-th harmonic number. This goes back to Moon in 1974 \cite{Moon1974}.

The path length of recursive trees was further investigated by Szymański \cite{Szymanski1990}, Dobrow and Smythe \cite{Dobrow1996}, Mahmoud \cite{Mahmoud1991} and Dobrow and Fill \cite{Dobrow1999}. Other types of linear preferential attachment trees have been considered as well: plane-oriented recursive trees \cite{Mahmoud1992}, binary search trees \cite{Regnier1989,Rosler1991}, and (weighted) $d$-ary increasing trees \cite{MunsoniusRuschendorf2011}. See also \cite{Sulzbach2017} for a general treatment of linear preferential attachment trees.

Linear preferential attachment trees can also be seen as varieties of increasing trees. Bergeron, Flajolet and Salvy \cite{Bergeron1992} showed that the path length grows in the order of $n\log n$ for rather general classes of increasing trees (that include those of interest in this paper as special cases). Let us finally also mention that the distribution of the path length has also been studied in
the very general setting of split trees \cite{BroutinHolmgren2012}.

By contrast, the external path length has received much less attention (see \cite{Javanian2010,Javanian2013} for two instances), but by comparing it to the ordinary path length we can now transfer results from one to the other.

In the mathematical literature, the Wiener index was introduced as the \emph{distance} of a graph in 1976 by Entringer, Jackson and Snyder~\cite{WienerIndexFirstMathEntringerJacksonSnyder1976}. See also the survey by Dobrynin, Entringer and Gutman \cite{WienerIndexTheoryAndApplicationsDobryninEntringerGutman}. For conditioned Bienaym\'e--Galton--Watson trees, the expectation of the Wiener index is of asymptotic order $n^{5/2}$, which was shown by Entringer, Meir, Moon and Sz\'ekely~\cite{Entringer1994}. After rescaling, the Wiener index converges to a functional of a Brownian excursion, first shown by Janson~\cite{Janson03}.

For trees of logarithmic height, the Wiener index is typically of order $n^2\log n$, and several classes of trees have been investigated. Neininger considered the Wiener index of recursive trees and binary search trees~\cite{Neininger}, Munsonius and Rüschendorf covered $d$-ary increasing trees~\cite{Munsonius}, and Munsonius studied the class of general split trees~\cite{Munsonius_split_trees}.

While there are a number of graph-theoretical results on the terminal Wiener index (see for instance the survey \cite{GutmanFurtulaSurvey2010}), it has not been given much attention for random trees. However, as we will see, results on the Wiener index can be carried over to the terminal Wiener index.

The paper is structured as follows: we start with a formal definition of the tree models and the quantities of interest and state the main results in the following section. An outline of the proofs, which are based on the analysis of suitable generating functions, is given in Section~\ref{sec:outline}, with a more detailed discussion of conditioned Bienaym\'e--Galton--Watson trees and linear preferential attachment trees in Sections~\ref{sec:BGW} and~\ref{sec:LPA} respectively. A brief section on possible extensions and generalisations concludes the paper.

\section{Definitions and main results}\label{sec:def_res}

We study two pairs of closely related distance-based parameters of rooted trees. The first pair consists of two depth parameters; the path length and the external path length.
The \textit{path length} (also known as total height) of a tree $T$ is defined as the sum of all distances between vertices in the tree and the root $r$, i.e.,
\[P(T)=\sum_{v\in V(T)}d_T(r,v).\]
The \textit{external path length} (terminal path length, total leaf height) is defined analogously, but with the summation restricted to leaves (vertices without children), i.e.,
\[P_E(T)=\sum_{l\in L(T)}d_T(r,l),\]
where $L(T)$ is the set of all leaves of the tree $T$. Normalizing these parameters by the number of vertices and the number of leaves yields the average depth of a uniformly chosen vertex and a uniformly chosen leaf, respectively.

We also consider the Wiener index and terminal Wiener index. The \textit{Wiener index} of a tree $T$ is defined as the sum of all pairwise distances between vertices in the tree, i.e.,
\[W(T)= \sum_{\{u,v\}\subseteq V(T)}d_T(u,v) = \frac{1}{2}\sum_{u,v\in V(T)}d_T(u,v).\]
Here, the first sum is over unordered pairs, the second over ordered pairs of vertices, which explains the factor $\frac12$.
The \textit{terminal Wiener index} (external Wiener index) restricts this summation to pairs of leaves, in analogy to the relation between the path length and the external path length, i.e.,
\[\wt= \sum_{\{u,v\}\subseteq L(T)}d_T(u,v) = \frac{1}{2}\sum_{u,v\in L(T)}d_T(u,v).\]
Normalizing these parameters by the number of summands (the number of pairs of vertices and pairs of leaves respectively) yields the average distance between two uniformly chosen vertices and two uniformly chosen leaves, respectively.

It is convenient, as in \cite{Janson03} and \cite{WAGNER20121502}, to introduce an auxiliary quantity relating the path length and the Wiener index, namely
\[Q(T):=\sum_{v,u\in V(T)}d_T(r,v\wedge u),\]
where $v\wedge u$ denotes the last common ancestor of $v$ and $u$. Using the tree structure one can rewrite the distance between two vertices $v$ and $u$ as
\begin{align*}
    d(v,u)&=d(v,v\wedge u)+d(v\wedge u,u)\\
    &=d(v,r)-d(v\wedge u,r)+d(u,r)-d(v\wedge u,r).
\end{align*}
For a tree $T$, summing this over all vertices yields 
\[W(T)=\frac{1}{2}\sum_{v,u\in V(T)}d(v,u)=|T|\left(\sum_{v\in V(T)}d(v,r)\right)-\left(\sum_{v,u\in V(T)}d(v\wedge u ,r)\right),\]
and so we get the relation
\begin{equation}\label{eq_W_relation}
    W(T)=|T|P(T)-Q(T).
\end{equation}
So in order to study $W(T)$ it is enough to study $P(T)$ and $Q(T)$.
Analogously we define 
\[\qt:=\sum_{l,k\in L(T)}d_T(r,l\wedge k),\]
and get the expression
\begin{equation}\label{eq_Wt_relation}
    \wt=\lt\pt-\qt.
\end{equation}

A \textbf{simply generated family of trees}, as introduced by Meir and Moon \cite{Meir1978}, is described by a weight sequence $\phi_0,\phi_1,\ldots$, which for each tree $T$ associates a weight
\begin{equation}\label{eq:tree-weight}
    \omega(T)=\prod_{v\in V(T)}\phi_{\deg_+(v)}=\prod_{i\geq 0}\phi_i^{c_i(T)},
\end{equation}
where $\deg_+(v)$ denotes the outdegree of a vertex $v$ and $c_i(T)$ is the number of vertices in $T$ with outdegree $i$. The generating function of these weights, $G(x)=\sum_T \omega(T)x^{|T|}$, then satisfies the equation
\begin{equation}
G(x)=x\Phi(G(x)),
    \label{g=xphig}
\end{equation}
where $\Phi(t)=\sum_{i\geq 0}\phi_i t^i$. Whenever the radius of convergence of $\Phi$ is positive, there exists a probability distribution $\xi$ such that the distribution on a family of simply generated trees of size $n$ with probabilities proportional to the weights is the same as the distribution of a random \textbf{conditioned Bienaym\'e--Galton--Watson tree} with offspring distribution $\xi$.  This type of random tree is constructed by a branching process 
that starts with a root vertex and then attaches a number of children to it according to $\xi$. This is then repeated independently for all successors. Conditioning this process on resulting in a finite tree with $n$ vertices yields a conditioned Bienaym\'e--Galton--Watson tree. 

Setting $\phi_i = \P(\xi = i)$ turns $\omega(T)$ into the probability that the (unconditioned) process yields the tree $T$. Translating from a general weight sequence to an equivalent probability distribution can be achieved by a transformation known as \emph{tilting}, where $\phi_i$ is replaced by $\phi_i \cdot ab^i$ for suitable constants $a$ and $b$---see for example Drmota's book \cite{Drmota2009} or Janson's survey \cite{Janson2012} for details.

In this paper, we will make the standard assumption that $\E\,\xi=1$, i.e., that the process is \emph{critical} (see again \cite{Drmota2009} or \cite{Janson2012}), and also assume that $\Var\,\xi=\sigma^2<\infty$. For technical convenience, we only consider the aperiodic case, where $\operatorname{span}(\xi) = \gcd \{ i: \P(\xi = i) > 0 \} = 1$. However, it would not be difficult to extend all results to the case where $\operatorname{span}(\xi) > 1$.

Under these conditions ($\E\,\xi=1$, $\Var\,\xi=\sigma^2<\infty$ and $\operatorname{span}(\xi) = 1$), the generating function $G(x)$ has a dominant square root singularity at $1$ and is otherwise analytic in what is known as a $\Delta$-domain. This is well known in the case where $\xi$ has exponential moments, but perhaps less well known under the weaker assumption that the variance is finite. The proof of the following result, which will be crucial for our approach, is due to Janson~\cite[Lemma A2]{Janson2006}.

\begin{lemma}\label{lem:bgw-gf}
    Consider conditioned Bienaym\'e--Galton--Watson trees with an offspring distribution $\xi$ that satisfies $\E\,\xi=1$, $\Var\,\xi=\sigma^2<\infty$ and $\operatorname{span}(\xi) = 1$, and let $G(x)$ be the associated generating function. Then $G(x)$ is analytic in a domain of the form
\begin{equation}\label{eq:delta-domain}
\{x \in \C \,:\, |x| < 1 + \epsilon, x \neq 1, |\operatorname{arg} (x-1)| > \delta \},
\end{equation}
    where $\epsilon > 0$ and $0 < \delta < \frac{\pi}{2}$, and it also satisfies
    \begin{equation*}
        G(x) = 1-\frac{\sqrt{2}}{\sigma}\sqrt{1-x}+o(|x-1|^{1/2})
    \end{equation*}
    as $x \to 1$ within this domain.
\end{lemma}

Let us now state the main results of this paper for Bienaym\'e--Galton--Watson trees.

\begin{theorem}\label{Theorem_simplygenerated}
    Let $T$ denote a conditioned Bienaym\'e--Galton--Watson tree of size $|T| = n$ with offspring distribution $\xi$ satisfying $\E\,\xi=1$, $\Var\,\xi=\sigma^2<\infty$ and $\operatorname{span}(\xi) = 1$. Set $\phi_0 = \P(\xi = 0)$. Then, as $n \to \infty$,
    \begin{align}
        \E\left(\frac{\pt}{\lt}-\frac{\p}{|T|}\right)&\to \frac{1}{\sigma^2},\label{path1Thm}\\
        \E\left(\left(\frac{\pt}{\lt}-\frac{\p}{|T|}\right)^2\right)&\to \frac{2\phi_0+(1-\phi_0)\sigma^2}{3\phi_0\sigma^4},\label{TheoremEqSecondMomentGW}
    \end{align}
    and 
    \begin{equation*}
        \frac{\pt}{\lt}-\frac{\p}{|T|}=O_p(1),    \end{equation*}
    where $O_p(1)$ denotes boundedness in probability.
\end{theorem}

\begin{theorem}\label{Theorem:BGW_Wiener}
    Let $T$ be as in Theorem~\ref{Theorem_simplygenerated}. Then, as $n\to\infty$,
        \begin{align}\label{eq:BGW_Wiener_Thm_eq1}
         \E\left(\frac{\wt}{\blt}-\frac{\w}{\bt}\right)&\to\frac{2}{\sigma^2},\\
         \E\left(\left(\frac{\wt}{\blt}-\frac{\w}{\bt}\right)^2\right)&\to\frac{4(14\phi_0+(1-\phi_0)\sigma^2)}{15\phi_0\sigma^4},\label{eq:BGW_Wiener_Thm_eq2}
    \end{align}
    and
    \begin{equation*}
        \frac{\wt}{\blt}-\frac{\w}{\bt}=O_p(1).
    \end{equation*}
\end{theorem}

\begin{remark}\label{remark:BGW_tree_is_path}
    When $T$ has only one leaf (so $T$ is a path rooted at one of the endpoints), $\wt=\blt=0$ and the quotient $\frac{\wt}{\blt}$ is not well-defined. Setting $\frac{\wt}{\blt}=0$ in this case resolves this issue since, by Lemma~\ref{lem:concentration-bgw}, the probability that it happens vanishes exponentially.
\end{remark}

\begin{remark}
    In principle, it would be possible to determine the asymptotic behaviour of higher moments as well. However, we were not able to identify general formulas for all moments that would potentially characterise the limiting distributions of the differences $\frac{\pt}{\lt}-\frac{\p}{|T|}$ or $\frac{\wt}{\blt}-\frac{\w}{\bt}$. This is left as an open problem.
\end{remark}

As a side result to the fact that the average distances to the root, i.e., $\frac{\p}{|T|}$ and $\frac{\pt}{\lt}$, are close to each other, we also obtain the limiting distribution of $\pt$. As mentioned in the introduction, $\frac{\p}{n^{3/2}}$ converges weakly to a multiple of the Airy distribution, which can be described in terms of the area under a Brownian excursion. We now find that this is also the case for $\frac{\pt}{n^{3/2}}$.

\begin{corollary}\label{cor:airy}
Let $T$ denote a conditioned Bienaym\'e--Galton--Watson tree of size $|T| = n$ with offspring distribution $\xi$ satisfying the same conditions as in Theorem~\ref{Theorem_simplygenerated}. 

Let the random variable $\eta$ be defined by
    \begin{equation*}
        \eta := 2 \int_0^1 e(t)dt,
    \end{equation*}
where $e(t)$, $0\leq t\leq 1$, is a normalized Brownian excursion. Then we have
\begin{equation*}
        \frac{\pt}{ n^{3/2}}\xrightarrow{d}\phi_0 \sigma^{-1} \eta.
\end{equation*}
\end{corollary}
The limiting distribution of $\frac{\w}{n^{5/2}}$ can, as mentioned in the introduction, also be described in terms of a Brownian excursion. This is true for $\frac{\wt}{n^{5/2}}$ as well.
\begin{corollary}\label{cor:BGW_Wiener}
    Let $T$ be as in Corollary~\ref{cor:airy}. Then
    \begin{equation*}
        \frac{\wt}{ n^{5/2}}\xrightarrow{d}\phi_0^2 \sigma^{-1} \zeta,
    \end{equation*}
    where $\zeta$ is the random variable
    \begin{equation*}
        \zeta := 2 \iint_{0<s<t<1}(e(s)+e(t)-2\min_{s\leq u\leq t} e(u))\,ds\,dt,
    \end{equation*}
    with $e(t)$ as in Corollary~\ref{cor:airy}.
\end{corollary}

\begin{remark}
It was shown by Marckert and Mokkadem~\cite{Marckert2003} that the height process of the leaves has the same limit as the height process associated with all vertices in the tree. In fact, this holds in greater generality for vertices with given outdegree. As such, it should not come as a surprise that the first order asymptotics of $\frac{\pt}{\lt}$ and $\frac{\p}{|T|}$ agree. However, it does not seem obvious that the difference is even bounded in probability.
\end{remark}

\textbf{Random linear preferential attachment trees} 
are trees constructed using a growth process, starting with a single root, where the probability of attachment to a vertex depends linearly on its number of children. For a real number $\alpha$, a random (linear) preferential attachment tree $\T_{\alpha,n}$ is obtained from $\T_{\alpha,n-1}$ by attaching the $n$-th vertex to $v\in V(\T_{\alpha,n-1})$ with probability
\begin{equation*}
    \frac{\alpha \deg_+(v)+1}{\sum_{u\in V(\T_{\alpha,n-1})}(\alpha \deg_+(u)+1)}=\frac{\alpha \deg_+(v)+1}{\alpha(n-1)+n},
\end{equation*}
where $\deg_+(v)$ is again the number of children of $v$.
To avoid negative probability of choosing a vertex $v$, this is restricted to $\alpha\in \{\ldots,-\frac{1}{3},-\frac{1}{2}\}\cup [0,\infty)$. This also excludes the case $\alpha=-1$, which simply produces a rooted path. In order to perform our analysis, we need to consider three separate cases: $\alpha=0$, $\alpha\in\{\ldots,-\frac{1}{3},-\frac{1}{2}\}$, and $\alpha>0$.
\begin{itemize}
    \item $\alpha=0$: 
    When $\alpha=0$, the choice between existing vertices is uniform and we get \textbf{random recursive trees}. They can be defined either using the growth process described above, or as labelled trees with vertices labelled $1,\ldots,n$ in an increasing fashion, chosen uniformly at random. These trees have been studied since at least 1967 (see Tapia and Myers~\cite{Tapia1967}). There is a bijection between trees of size $n$ and permutations of $n-1$ elements, and thus the number of recursive trees of size $n$ is $(n-1)!$.
\item $\alpha\in\{\ldots,-\frac{1}{3},-\frac{1}{2}\}$:
When $\alpha=-\frac{1}{d}$ for a positive integer $d$, we get \textbf{random $d$-ary increasing trees}. From a more combinatorial viewpoint, these can be seen as trees where each vertex has $d$ places where a child can be attached (thus up to $d$ children), and vertices are again labelled $1,\ldots,n$ in an increasing fashion.

\item $\alpha>0$: these are preferential attachment trees in the sense that vertices with many children have a greater probability of attracting new vertices. These trees are also known as \textbf{random generalised plane-oriented recursive trees} (GPORTs), as they extend the special case of \textbf{plane oriented recursive trees} (introduced by Szymański \cite{Szymanski1987}), or PORTs, which corresponds to $\alpha=1$. 
\end{itemize}

All three cases can also be obtained from a combinatorial construction similar to that of simply generated families of trees: on the set of all rooted ordered trees with vertices labelled $1,\ldots,n$ in an increasing fashion (i.e., the labels increase as one moves away from the root, see also~\cite{Bergeron1992}), one defines a weight function as in~\eqref{eq:tree-weight}, with $\phi_i$ as follows:
\begin{itemize}
\item $\phi_i = \frac{1}{i!}$ to obtain recursive trees,
\item $\phi_i = \binom{d}{i}$ to obtain $d$-ary increasing trees,
\item $\phi_i = \binom{i - 1 + \frac{1}{\alpha}}{i}$ to obtain GPORTs.
\end{itemize}
Random trees in each of these cases are obtained by choosing trees of a given size with probabilities proportional to the weights. As in the case of simply generated families, one could multiply $\phi_i$ by $ab^i$ for arbitrary positive constants $a$ and $b$ to obtain a probabilistically equivalent choice of weights. While the choices of weights might look rather different, we remark that they are always equal to $\frac{1}{i!}\prod_{j=0}^{i-1} (1+\alpha j)$ times a factor of the form $b^i$.

It was shown by Panholzer and Prodinger~\cite{PanholzerProdinger2007} that these are the only choices of weights for which the probabilistic model is equivalent to a suitable growth process. 

Aiming towards a similar analysis as in the simply generated case, we consider exponential generating functions of the form
\[\sum_T \omega(T) \frac{x^{|T|}}{|T|!},\]
where the sum is over all rooted ordered increasing trees. 
We use $\Gr$ for the exponential generating function associated with recursive trees, $\Gd$ for $d$-ary increasing trees and $\Gp$ for GPORTs. Each of these satisfies a first order differential equation (see \cite{Drmota2009}).
Specifically,
\begin{equation}
\Gr'(x)=e^{\Gr(x)}\label{recursiveODE}
\end{equation}
for recursive trees, with the explicit solution
\[\Gr(x)=\log\left(\frac{1}{1-x}\right).\]
Likewise,
\begin{equation}
    \Gd'(x)=(1+\Gd(x))^d \label{daryODE}
\end{equation}
with the explicit solution
\begin{equation*}
    \Gd(x)=(1-(d-1)x)^{\frac{1}{1-d}}-1.
\end{equation*}
Finally, 
\begin{equation}
    \Gp'(x)=(1-\Gp(x))^{-\frac{1}{\alpha}} \label{prefattODE}
\end{equation}
with the explicit solution
\begin{equation*}
    \Gp(x)=1-\left(1-\left(1+\frac{1}{\alpha}\right)x\right)^{\frac{\alpha}{1+\alpha}}.
\end{equation*}
These differential equations and the slightly different shapes of their solutions are a main reason for differentiating between the three  cases.

Let us now state the main result for linear preferential attachment trees, which parallels Theorem~\ref{Theorem_simplygenerated} and Theorem~\ref{Theorem:BGW_Wiener}. 

\begin{theorem}\label{sats_alltogether}
Let $T$ denote a random linear preferential attachment tree $\T_{\alpha,n}$ of size $|T| = n$, with $\alpha \in \{\ldots,-\frac{1}{3},-\frac{1}{2}\}\cup [0,\infty)$.
 Then, as $n \to \infty$,
    \begin{align*}
        \E\left(\frac{\pt}{\lt}-\frac{\p}{|T|}\right)&\to \frac{1}{(1+\alpha)(2+\alpha)},\\
        \E\left(\left(\frac{\pt}{\lt}-\frac{\p}{|T|}\right)^2\right)&\to \frac{1}{(1+\alpha)^2(2+\alpha)^2},
    \end{align*}
    and
   \begin{equation*}
       \frac{\pt}{\lt}-\frac{\p}{|T|}\xrightarrow{P}\frac{1}{(1+\alpha)(2+\alpha)},
   \end{equation*} 
   where $\xrightarrow{P}$ denotes convergence in probability.
\end{theorem}
\begin{theorem}\label{Theorem:increasing_Wiener}
Let $T$ be as in Theorem~\ref{sats_alltogether}.
 Then, as $n \to \infty$,
    \begin{align*}
        \E\left(\frac{\wt}{\blt}-\frac{\w}{\bt}\right)&\to \frac{2}{(1+\alpha)(2+\alpha)},\\
        \E\left(\left(\frac{\wt}{\blt}-\frac{\w}{\bt}\right)^2\right)&\to \frac{4}{(1+\alpha)^2(2+\alpha)^2},
    \end{align*}
    and
   \begin{equation*}
       \frac{\wt}{\blt}-\frac{\w}{\bt}\xrightarrow{P}\frac{2}{(1+\alpha)(2+\alpha)}.
   \end{equation*} 
\end{theorem}

\begin{remark}\label{remark:lpat_tree_is_path}
As before, we set $\frac{\wt}{\blt}=0$ if numerator and denominator are both $0$. By Lemma~\ref{lem:concentration-bgw}, this only happens with exponentially small probability.
\end{remark}
An analogue of Corollary~\ref{cor:airy} follows as well. For random linear preferential attachment trees, $\frac{\p - \E(\p)}{n}$ converges in distribution to a limiting random variable $X_{\alpha}$ (whose distribution depends on $\alpha$) as $n \to \infty$ (in fact, even almost sure convergence holds, see \cite{Sulzbach2017}). 

\begin{corollary}\label{cor:Xalpha}
Let $T$ denote a random linear preferential attachment tree $\T_{\alpha,n}$ of size $|T| = n$, with $\alpha \in \{\ldots,-\frac{1}{3},-\frac{1}{2}\}\cup [0,\infty)$, and let $X_{\alpha}$ be as described above.  Then, as $n \to \infty$,
\[\frac{\pt - \E(\pt)}{n} \xrightarrow{d} \frac{1+\alpha}{2+\alpha} X_{\alpha}.\]
\end{corollary}
Provided that a similar limit law holds for the Wiener index, that is, $\frac{\w-\E\left(\w\right)}{n^2}$ converges to a random variable $Y_\alpha$ (with a distribution depending on $\alpha$), a limit law for the terminal Wiener index follows as well. This was shown by Neininger~\cite{Neininger} in the special cases of recursive trees and binary increasing trees, and the more general cases $\alpha \in \{\ldots,-\frac{1}{3},-\frac{1}{2}\}\cup \{0,1,2,3,\ldots\}$ were treated by Munsonius and R\"uschendorf~\cite{MunsoniusRuschendorf2011}. An analogous limit theorem most probably also holds for positive non-integer values of $\alpha$, but it does not seem to be mentioned anywhere in the literature.
\begin{corollary}\label{cor:increasing_Wiener}
    Let $T$ be as in Corollary~\ref{cor:Xalpha}, and let $Y_{\alpha}$ be as described above. Then, as $n\to\infty$,
    \[\frac{\wt - \E(\wt)}{n^2} \xrightarrow{d} \left(\frac{1+\alpha}{2+\alpha}\right)^2 Y_{\alpha}.\]
\end{corollary}

\begin{remark}
    In~\cite{TerminalWienerIndexGutman}, the terminal Wiener index is not defined using leaves, since it is defined for general graphs, but rather using vertices of degree 1. In the context of rooted trees, the definition used here seems more natural, and asymptotically, the difference is actually negligible, as will be explained in the following. When the root has more than one child, the definitions are equivalent, the only difference arises when the root has degree $1$.

    Let $W^*_E(T)$ and $L^*(T)$ denote the terminal Wiener index and the set of leaves respectively, both including the root $r$ when its outdegree is $1$. The contribution of the root is exactly the external path length times the indicator that it has degree 1. Let us write $\mathbf{1}_r$ for this indicator, so that
    \begin{align*}
        |L^*(T)|&=\lt+\mathbf{1}_r,\\
        W^*_E(T)&=\wt+\mathbf{1}_r\pt.
    \end{align*}
    Whenever $\lt\geq 2$, we can consider the difference of the two averages and obtain the upper bound
    \begin{align*}
        \bigg|\frac{W^*_E(T)}{\binom{|L^*(T)|}{2}}-\frac{\wt}{\blt}\bigg|&=\bigg|\wt\left(\frac{\blt-\binom{|L^*(T)|}{2}}{\blt\binom{|L^*(T)|}{2}}\right)+\frac{\mathbf{1}_r\pt}{\binom{|L^*(T)|}{2}}\bigg|\\
        &\leq \frac{4\wt}{\lt(\lt^2-1)}+\frac{2\pt}{\lt(\lt-1)},
    \end{align*}
    by the triangle inequality and the fact that \[\bigg|\frac{\blt-\binom{|L^*(T)|}{2}}{\blt\binom{|L^*(T)|}{2}}\bigg|\leq \bigg|\frac{\blt-\binom{\lt+1}{2}}{\blt \binom{\lt+1}{2}}\bigg| = \frac{4}{\lt(\lt^2-1)}.\] 
    As will also be used later, $\frac{\lt}{n}$ converges in probability to a constant (in all cases that we consider). At the same time, $\frac{\wt}{n^3}$ and $\frac{\pt}{n^2}$ both tend to $0$ in probability. This can be seen for example by a basic application of the Markov inequality, making use of the fact that the expected value of the Wiener index is of order $n^{5/2}$ (for conditioned Bienaym\'e--Galton--Watson trees) and $n^2 \log n$ (for linear preferential attachment trees) respectively, while the expected value of the path length is of order $n^{3/2}$ (for conditioned Bienaym\'e--Galton--Watson trees) and $n \log n$ (for linear preferential attachment trees) respectively.

    As a consequence, both $\frac{4\wt}{\lt(\lt^2-1)}$ and $\frac{2\pt}{\lt(\lt-1)}$ go to $0$ in probability. So the same is true for the difference of the two averages, which is thus negligible.
\end{remark}

\section{Proof outline and preliminaries}\label{sec:outline}
The proofs of Theorems \ref{Theorem_simplygenerated} and \ref{sats_alltogether} are based on the analysis of multivariate generating functions encoding the relevant statistics $|T|$, $|L(T)|$, $\p$, $\pt$, $\q$ and $\qt$ (together giving $\w$ and $\wt$).

The functional and differential equations satisfied by the various generating functions reflect the recursive structure of trees: a tree can be decomposed into a root and smaller trees (branches) that are attached to it. 

In order to study the path length, external path length, Wiener index and terminal Wiener index we will therefore consider the following recursions that also make use of this branch structure. If $T$ has more than one vertex and $T_1,\,T_2,\ldots,T_m$ are the branches of $T$ (the rooted trees that remain as components when the root is removed), then we have the recursions
\begin{align}\label{eq:leafrec}
    |L(T)| &= \sum_{i=1}^m |L(T_i)|, \\
\p &=\sum_{i=1}^m\left(P(T_i)+|T_i|\right) \label{eq:plrec}
\end{align}
and 
\begin{equation}\label{eq:eplrec}
\pt=\sum_{i=1}^m\left(P_E(T_i)+|L(T_i)|\right).
\end{equation}
Instead of obtaining similar recursions for the Wiener index, we do so for $Q(T)$ and $\qt$. If $T_1,\,T_2,\ldots,T_m$ are the branches of $T$, then we have the recursion
\begin{equation}
    Q(T)=\sum_{i=1}^m\left(Q(T_i)+|T_i|^2\right).
    \end{equation}
The external version is similarly
\begin{equation}
\qt=\sum_{i=1}^m\left(\qt(T_i)+|L(T_i)|^2\right).
\end{equation}

In the Bienaym\'e--Galton--Watson case, this leads to an implicit functional equation extending \eqref{g=xphig} to the multivariate setting. For the different varieties of increasing trees, we work with an exponential multivariate generating function and obtain partial differential equations that extend the ODEs (\ref{recursiveODE}), (\ref{daryODE}) and (\ref{prefattODE}).
Repeated differentiation then allows us to determine generating functions for the first and second moments of the random variables
\[|T|\pt-\lt P(T)\]
and
\[\bt\wt-\blt \w .\]
Using singularity analysis and the known behaviour of the univariate generating functions for each tree family, we are able to determine the asymptotic behaviour of these moments. The calculations related to the generating functions are carried out in SageMath \cite{sagemath}; the complete code used for this project is available in a public repository.\footnote{All code used for symbolic computations is available at:
\url{https://github.com/JacobLundblad/ExternalPathLengthAndTerminalWienerIndex}.}
Using concentration results for the number of leaves $\lt$ allows for these moment asymptotics to be used to determine the asymptotic behaviour of the first and second moments of the quotients
\[\frac{|T|\pt-\lt P(T)}{|T|\lt}\]
and
\[\frac{\bt \wt-\blt\w}{\bt\blt}.\]
To handle these ratios, we rely on the technical lemma stated below.
\begin{lemma}
\label{XoverYLemma}
    Let $X_n$ and $Y_n$ be random variables. 
    Assume that $\frac{\E(X_n)}{\mu_n}$ and $\frac{\E(X_n)}{\E(|X_n|)}$ are bounded below by positive constants, that $|X_n|\leq h(n)$ for some function $h$, and that $|Y_n|\geq 1$. If $Y_n$ satisfies
    \[\P\left(|Y_n-\mu_n|>g(n)\right) \leq f(n)\] 
    for some functions $f$ and $g$ such that $f(n)=o(1)$, $h(n)=o(1/f(n))$ and $g(n)=o(\mu_n)$, then
    \begin{equation*}
        \E\left(\frac{X_n}{Y_n}\right)\sim \frac{\E(X_n)}{\mu_n}\quad\text{as }n\to\infty.
    \end{equation*}
\end{lemma}
\begin{remark}
    Note that we do not require $\E(Y_n)=\mu_n$.
\end{remark}
\begin{proof}
    Consider the event $A=\{|Y_n-\mu_n| > g(n)\}$, and split the expectation as follows:
    \begin{equation}\label{eq:A-split}
        \E\left(\frac{X_n}{Y_n}\right)=\E\left(\frac{X_n}{Y_n}\,\Big|\,A\right)\P(A)+\E\left(\frac{X_n}{Y_n}\,\Big|\,\Bar{A}\right)\P(\Bar{A}).
    \end{equation}
    Starting with the first term, we can use the bounds for $X_n$ and $Y_n$ to obtain
    \begin{align*}
        \E\left(\big|\frac{X_n}{Y_n}\big| \,\Big|\,A\right)\leq\E\left(|X_n|\big|A\right)\leq h(n).
    \end{align*}
    Thus, we have  
    \begin{equation*}
        \Big|\E\left(\frac{X_n}{Y_n} \,\Big|\,A\right)\P(A)\Big| \leq \E\left(\big|\frac{X_n}{Y_n}\big| \,\Big|\,A\right)\cdot \P(A) \leq h(n)\cdot f(n) =o(1).
    \end{equation*}
    For the second term, we note that the expectation of $X_n$ can be written as
    \[\E(X_n)=\E(X_n|A)\P(A)+\E(X_n|\Bar{A})\P(\Bar{A}),\]
    and so
    \[\E(X_n|\Bar{A})=\frac{\E(X_n)-\E(X_n|A)\P(A)}{\P(\Bar{A})}.\]
    Using the bounds on $X_n$ and $\P(A)$, we get that
    \[\E(X_n|\Bar{A})=\frac{\E(X_n)}{\P(\Bar{A})}+O\left(\frac{h(n)\cdot f(n)}{1-f(n)}\right)\sim \E(X_n),\]
    since $\P(\Bar{A})\xrightarrow{n\to\infty}1$ by the assumption that $f(n)=o(1)$. Likewise, we have
    \[\E\big(|X_n| \, \big|\, \Bar{A}\big) \sim \E(|X_n|).\] 
    Under the condition $\Bar{A}$, we have $Y_n=\mu_n+O(g(n))$. Thus, it follows that
    \begin{align*}
        \E\left(\frac{X_n}{Y_n} \,\Big|\,\Bar{A}\right)&=\E\left(\frac{X_n}{\mu_n+O(g(n))} \,\Big|\,\Bar{A}\right)\\
        &=\E\left(\frac{X_n}{\mu_n} + O \Big( \frac{|X_n|g(n)}{\mu_n^2} \Big) \,\Big|\,\Bar{A}\right)\\
        &=\frac{\E(X_n  \,|\,\Bar{A})}{\mu_n} + O \Big( \frac{\E\big(|X_n| \, \big|\, \Bar{A}\big) g(n)}{\mu_n^2} \Big) \\ 
        &\sim \frac{\E(X_n)}{\mu_n}.
    \end{align*}
    Combining \eqref{eq:A-split} with our estimates for $\E\left(\frac{X_n}{Y_n}\,\Big|\,A\right)$ and $\E\left(\frac{X_n}{Y_n}\,\Big|\,\Bar{A}\right)$ finally gives us
    \[\E\left(\frac{X_n}{Y_n}\right) \sim \frac{\E(X_n)}{\mu_n},\]
    which proves the lemma.
\end{proof}

For some calculations we will also use the following lemma.
\begin{lemma}
    \label{Lemma_highermomentsconsentration}
    Let $Y_n$ be a non-negative random variable satisfying 
    \[\P(|Y_n-\mu_n|>g(n)) \leq f(n)\]
    for some functions $\mu_n$, $f$ and $g$, and assume that $Y_n\leq h(n)$ and $0\leq\mu_n\leq h(n)$ for some function $h$. Then, 
    \begin{equation*}
        \P\left(|Y_n^2-\mu_n^2|>2g(n)h(n)\right) \leq f(n).
    \end{equation*}
\end{lemma}

\begin{proof}
We have    
    \begin{align*}
        f(n) &\geq\P(|Y_n-\mu_n|>g(n)) \\
        &=\P(2h(n)|Y_n-\mu_n|>2g(n)h(n))\\
        &\geq \P\left((Y_n+\mu_n)|Y_n-\mu_n|>2g(n)h(n)\right) \\
        &=\P\left(|Y_n^2-\mu_n^2|>2g(n)h(n)\right).\qedhere
    \end{align*}
\end{proof}
In order to use Lemma~\ref{XoverYLemma} and Lemma~\ref{Lemma_highermomentsconsentration}, we will need the following concentration results for the number of leaves for the different types of trees.

\begin{lemma}\label{lem:concentration-bgw}

Let $T$ denote a conditioned Bienaym\'e--Galton--Watson tree of size $|T| = n$ with offspring distribution $\xi$ satisfying $\E\,\xi=1$, $\Var\,\xi=\sigma^2<\infty$ and $\operatorname{span}(\xi) = 1$. The number of leaves $\lt$ satisfies the concentration inequality
\begin{equation*}
 \P(|\lt-\phi_0n| > x)\leq \kappa \sqrt{n} e^{-\lambda x^2/n}
\end{equation*}
for certain positive constants $\kappa$ and $\lambda$ that only depend on $\xi$.
\end{lemma}

\begin{proof}
We take the viewpoint of the Bienaym\'e--Galton--Watson process. Let $\xi_i$ be the offspring of the $i$-th vertex, so that $\xi_1,\xi_2,\ldots$ are independent random variables with the same distribution $\xi$. Conditioning on the event that $S_n=\sum_{i=1}^n \xi_i = n-1$ while $S_k \geq k$ for $k < n$, we get a conditioned Bienaym\'e--Galton--Watson tree of size $n$. 
Note that the leaves are precisely the vertices for which $\xi_i = 0$. By the cycle lemma, we can equivalently condition on $S_n=\sum_{i=1}^n \xi_i = n-1$ only, since we are only interested in the number of $\xi_i$ that are equal to $0$.

Define $\hat{\xi_i}$ as $\xi_i$ conditioned on being non-zero, and set $\Hat{S}_{n-l}=\sum_{i=1}^{n-l}\hat{\xi}_i$. Then we have 
\begin{align*}
    \P(\lt=l)&=\P(\text{exactly $l$ of the $\xi_i$ are $0$}\,\big|\,S_n=n-1)\\
    &=\frac{\P(S_n=n-1,\text{exactly $l$ of the $\xi_i$ are $0$})}{\P(S_n=n-1)}\\
    &=\binom{n}{l}\frac{\P(S_{n-l}=n-1, \xi_i\neq 0\text{ for }i\leq n-l,\,\xi_i= 0\text{ for }i>n-l)}{\P({S_n}=n-1)}\\
    &=\binom{n}{l}\phi_0^l\frac{\P(S_{n-l}=n-1, \xi_i\neq 0\text{ for }i\leq n-l)}{\P({S_n}=n-1)}\\
    &=\binom{n}{l}\phi_0^l(1-\phi_0)^{n-l}\frac{\P({S_{n-l}}=n-1\,\big|\,\xi_i\neq 0 \text{ for }i\leq n-l)}{\P(S_n=n-1)}\\
    &=\binom{n}{l}\phi_0^l(1-\phi_0)^{n-l}\frac{\P(\hat{S}_{n-l}=n-1)}{\P(S_n=n-1)}\\
    &\leq \P(\text{Bi}(n,\phi_0)=l)\cdot c\sqrt{n}
\end{align*}
for some constant $c$, where $\text{Bi}(n,\phi)$ denotes the binomial distribution with parameters $n$ and $\phi$. The last inequality bounds the numerator by 1 and uses the classical local central limit theorem for discrete distributions, see e.g.~\cite[Chapter 9]{Gnedenko1968}. 
Thus we get
\begin{align*}
    \P(|\lt-\phi_0n| > x)&=\sum_{l=1}^{n-1}\P(\lt=l)I_{|l-\phi_0n| > x}\\
    &\leq c\sqrt{n} \sum_{l=1}^{n-1}\P(\text{Bi}(n,\phi_0)=l)I_{|l-\phi_0n|> x}\\
    &=c\sqrt{n} \,\P(|\text{Bi}(n,\phi_0)-\phi_0n|> x).
\end{align*}
Combining this with the Chernoff bound
\[\P(|\text{Bi}(n,\phi_0)-\phi_0n|\geq\delta\phi_0n)\leq2e^{-\phi_0n\frac{\delta^2}{3}},\]
which is valid for any $\delta\in(0,1)$, 
yields the desired concentration inequality, proving the lemma.
\end{proof}

We obtain a similar lemma for the random linear preferential attachment trees.
In the special case of recursive trees, such a concentration result was proven by Zhang in~\cite{Zhang2015}.

\begin{lemma}\label{lem:concentration-lpat}
Let $T$ be a random linear preferential attachment tree $\mathcal{T}_{\alpha,n}$ of size $|T| = n$ with $\alpha \in \{\ldots, -\frac13, -\frac12 \} \cup [0,\infty)$. The number of leaves $\lt$ satisfies the concentration inequality
\begin{equation*}
    \P \Big( \Big|\lt - \frac{1+\alpha}{2+\alpha} n \Big| > x \Big) \leq \kappa e^{-\lambda x^2/n}
\end{equation*}
for certain positive constants $\kappa$ and $\lambda$ that only depend on $\alpha$.   
\end{lemma}

\begin{proof}
Let us consider the growth process that attaches a new leaf in each step as described in Section~\ref{sec:def_res}. Suppose that $L_k$ is the number of leaves after $k$ steps. Then we have (deterministically) $L_1 = L_2 = 1$, and the random variable $\lt$ we are interested in is simply $L_n$. After $k$ steps, there are $L_k$ leaves and $k-L_k$ internal vertices. By definition of the process, the probability to attach to a fixed leaf in the $(k+1)$-th step is
$\frac{1}{\alpha (k-1) + k}$. Thus, the $(k+1)$-th vertex is attached to a leaf with probability $\frac{L_k}{\alpha(k-1) + k}$, and to an internal vertex 
with probability $1 - \frac{L_k}{\alpha(k-1) + k}$. In the former case, the number of leaves stays the same, in the latter case it increases by $1$. Hence we have
\begin{equation*}
    \E \big( L_{k+1} \,\big|\, L_k \big) = L_k + 1 - \frac{L_k}{\alpha(k-1) + k} = \frac{(1+\alpha)(k-1)}{\alpha(k-1)+k} L_k + 1 = \frac{k-1}{k- \frac{\alpha}{1+\alpha}} L_k + 1.
\end{equation*}
We multiply by $\frac{\Gamma(k+1/(1+\alpha))}{\Gamma(k)} = \frac{\Gamma(k+1-\alpha/(1+\alpha))}{\Gamma(k)}$ to find
\begin{equation*}
\E \Big( \frac{\Gamma(k+1/(1+\alpha))}{\Gamma(k)} L_{k+1} \,\big|\, L_k \Big) = \frac{\Gamma(k-\alpha/(1+\alpha))}{\Gamma(k-1)} L_{k} + \frac{\Gamma(k+1/(1+\alpha))}{\Gamma(k)}.
\end{equation*}
Thus, if we set
\begin{equation*}
    X_k = \frac{\Gamma(k-\alpha/(1+\alpha))}{\Gamma(k-1)} L_{k} - \sum_{j=1}^{k-1} \frac{\Gamma(j+1/(1+\alpha))}{\Gamma(j)}
\end{equation*}
for $k \geq 2$, then $X_k$ is a martingale with $X_2 = 0$. Note also that
\begin{equation*}
X_{k+1} - X_k = \begin{cases} \frac{\Gamma(k-\alpha/(1+\alpha)) L_k}{(1 + \alpha) \Gamma(k)} & L_{k+1} = L_k + 1, \\ \frac{\Gamma(k-\alpha/(1+\alpha)) (L_k - (1 + \alpha) k + \alpha)}{(1 + \alpha) \Gamma(k)} & L_{k+1} = L_k, \end{cases}
\end{equation*}
thus (as $L_k \leq k$) we have $|X_{k+1} - X_k| = O(k^{1/(1+\alpha)})$ by Stirling's formula. In other words, there exists a constant $C$ such that $|X_{k+1} - X_k| \leq Ck^{1/(1+\alpha)}$ for all $k \geq 2$. By the Azuma--Hoeffding inequality, we have, with
\begin{equation*}
d_n^2 = \sum_{j=2}^{n-1} C^2 j^{2/(1+\alpha)} = O \big( n^{1+2/(1+\alpha)} \big),\end{equation*}
that
\begin{equation*}
\P ( |X_n| > t ) \leq 2e^{-t^2/(2d_n^2)}
\end{equation*}
for all $t \geq 0$. Let us write $a_n = \frac{\Gamma(n-\alpha/(1+\alpha))}{\Gamma(n-1)}$ and $
b_n = \sum_{j=1}^{n-1} \frac{\Gamma(j+1/(1+\alpha))}{\Gamma(j)}$, so that $X_n = a_n L_n - b_n$. Note that (again by Stirling's formula)
\begin{equation*}
    a_n = \frac{\Gamma(n-\alpha/(1+\alpha))}{\Gamma(n-1)} =  n^{1/(1+\alpha)}( 1 + O(n^{-1}))
\end{equation*}
and thus
\begin{equation*}
b_n = \sum_{j=1}^{n-1} \frac{\Gamma(j+1/(1+\alpha))}{\Gamma(j)} = \frac{1+\alpha}{2+\alpha} n^{(2+\alpha)/(1+\alpha)} + O\big(n^{1/(1+\alpha)}\big) =
\frac{1+\alpha}{2+\alpha} a_n (n + O(1)).
\end{equation*}
Putting everything together, we obtain the estimate
\begin{align*}
\P \Big( \Big|L_n - \frac{1+\alpha}{2+\alpha} n \Big| > x \Big)
&= \P \Big( \Big|a_n L_n - \frac{1+\alpha}{2+\alpha} na_n  \Big| > a_n x \Big) \\
&\leq \P \Big( \Big|a_n L_n - b_n \Big| > a_n x - \Big|b_n - \frac{1+\alpha}{2+\alpha} na_n\Big| \Big) \\
&= \P \Big( |X_n| > a_n x - \Big|b_n - \frac{1+\alpha}{2+\alpha} na_n\Big| \Big) \\
&\leq 2 \exp \Big( - \frac{(a_n x - |b_n - \frac{1+\alpha}{2+\alpha} na_n|)^2}{2d_n^2} \Big) \\
&= 2 \exp \Big( - \frac{a_n^2 ( x - O(1))^2}{2d_n^2} \Big),
\end{align*}
which is valid at least for $x > \Big|\frac{b_n}{a_n} - \frac{1+\alpha}{2+\alpha} n\Big| = O(1)$. Since $\frac{d_n^2}{a_n^2} = O(n)$, it follows that there exist constants $\kappa, \lambda > 0$ such that
\begin{equation*}
    \P \Big( \Big|L_n - \frac{1+\alpha}{2+\alpha} n \Big| > x \Big) \leq \kappa e^{-\lambda x^2/n}
\end{equation*}
holds for all $n$ and all $x > 0$, which proves the lemma.
\end{proof}

\section{Bienaym\'e--Galton--Watson trees}\label{sec:BGW}

In this section, we discuss details of the proof of our main results in the case of conditioned Bienaym\'e--Galton--Watson trees. The main ingredient for the proofs of Theorem~\ref{Theorem_simplygenerated} and Theorem~\ref{Theorem:BGW_Wiener} is the following.
\begin{lemma}\label{Lemma_simplygenerated}
Let $T$ be as in Theorem~\ref{Theorem_simplygenerated}. Then we have 
\begin{equation*}
    \E\left(|T|\pt-\lt\p\right) \sim \frac{\phi_0}{\sigma^2}n^2
\end{equation*}    
and
\begin{equation*}
   \E\left(\left(|T|\pt-\lt\p\right)^2\right) \sim \frac{\phi_0(2\phi_0+(1-\phi_0)\sigma^2)}{3\sigma^4}n^4
\end{equation*}
as $n \to \infty$.
Similarly,
\begin{equation*}
    \E\left(\bt\wt-\blt\w\right)\sim\frac{\phi_0^2}{2\sigma^2}n^4
\end{equation*}
and
\begin{equation*}
        \E\left(\left(\bt\wt-\blt\w\right)^2\right)\sim\frac{\phi_0^3(14\phi_0+(1-\phi_0)\sigma^2)}{60\sigma^4}n^8
\end{equation*}
as $n\to\infty$.

\end{lemma}

\begin{proof}
Aiming to extend \eqref{g=xphig} to include all quantities of interest, namely $|T|$, $\lt$, $\p$, $\pt$, $\q$, $\qt$, as well as $|T|^2$ and $\lt^2$, we consider the generating function

\[F(x,y,z,u,v,w,s,t)=\sum_T \omega(T) x^{|T|}y^{\lt} z^{\pt} u^{\p}v^{\qt}w^{\q}s^{|T|^2}t^{\lt^2}.\]
This allows us to use the recursions~\eqref{eq:leafrec},~\eqref{eq:plrec} and~\eqref{eq:eplrec} for $\lt$, $\p$ and $\pt$ to get the following functional equation:
\begin{equation}    F(x,y,z,u,v,w,1,1)=x\Phi\left(F(xu,yz,z,u,v,w,w,v)\right)+\phi_0 x(y-1).\label{functionalExpression}
\end{equation}
Note here that the final term $\phi_0 x(y-1)$ accounts for the fact that the single vertex tree has one leaf. 

Let us now introduce differential operators $\delta_x$, $\delta_y$, $\delta_z$, $\delta_u$, $\delta_v$, $\delta_w$, $\delta_s$ and $\delta_t$ as follows: for a function $A(x,y,z,u,v,w,s,t)$, we define
\[(\delta_x A)(x,y,z,u,v,w,s,t)=\left(x\cdot\frac{\partial}{\partial x}A\right)(x,y,z,u,v,w,s,t).\]
The operators $\delta_y$, $\delta_z$, $\delta_u$, $\delta_v$, $\delta_w$, $\delta_s$ and $\delta_t$ are defined analogously. Note that these operators commute. For a polynomial $p=\sum_{i=0}^k c_ix^{x_i}y^{y_i}z^{z_i}u^{u_i}v^{v_i}w^{w_i}s^{s_i}t^{t_i}$, where $c_i$ are arbitrary real coefficients and $x_i,y_i,z_i,u_i,v_i,w_i,s_i,t_i$ are non-negative integers for all $i$, we define a new operator $\Delta_p$, where $p$ acts as an indicator of which differential operators to apply. We define 
\[(\Delta_pA)(x,y,z,u,v,w,s,t)=\left(\sum_{i=0}^kc_i\delta_x^{x_i}\delta_y^{y_i}\delta_z^{z_i}\delta_u^{u_i}\delta_v^{v_i}\delta_w^{w_i}\delta_s^{s_i}\delta_t^{t_i} A\right)(x,y,z,u,v,w,s,t).\]
Now, we can express the generating functions for the expressions in our lemma by applying $\Delta_p$ to $F$ for a suitable $p$. In particular, for $p_1=xz-yu$ and $p_2=(xz-yu)^2$,
\[
    (\Delta_{p_1}F)(x,1,1,1,1,1,1,1) =\sum_T \omega(T) \big(|T|\pt-\lt\p\big) x^{|T|}
\]
and
\[
    (\Delta_{p_2}F)(x,1,1,1,1,1,1,1) =\sum_T \omega(T) \big(|T|\pt-\lt\p\big)^2 x^{|T|}.
\]
Motivated by the relations~\eqref{eq_W_relation} and~\eqref{eq_Wt_relation}, define also $p_3=\binom{x}{2}(yz-v)-\binom{y}{2}(xu-w)$ and $p_4=\left(\binom{x}{2}(yz-v)-\binom{y}{2}(xu-w)\right)^2$. We then have 
\[  (\Delta_{p_3}F)(x,1,1,1,1,1,1,1)=\sum_T \omega(T)\left(\binom{|T|}{2}\wt-\binom{\lt}{2}\w\right)x^{|T|}\]
and
\[  (\Delta_{p_4}F)(x,1,1,1,1,1,1,1)=\sum_T \omega(T)\left(\binom{|T|}{2}\wt-\binom{\lt}{2}\w\right)^2x^{|T|}.\]

This also motivates the definition of $\Delta_p$, as we note that writing functions like $\Delta_{p_4} F$ in terms of the various $\delta$-operators would be somewhat cumbersome. To then get the expectations in the lemma we need to divide the coefficient of $x^n$ in these generating functions by the coefficient of $x^n$ in
\[G(x):= F(x,1,1,1,1,1,1,1)=\sum_T \omega(T) x^{|T|}.\]
As previously stated, $G(x)=x\Phi(G(x))$ and thus
\[\Phi\big(F(xu,yz,z,u,v,w,w,v)\big)\big|_{y=z=u=v=w=1}=\Phi(G(x))=\frac{G(x)}{x}.\]
Applying $\Delta_{p_i}$ to the right-hand side of \eqref{functionalExpression} allows us to express all generating functions in terms of $G$ and its derivatives.

For $(\Delta_{p_1}F)(x,1,1,1,1,1,1,1)$, it is possible to perform the calculations by hand, as illustrated in Example~\ref{example_GW}. For the rest we use SageMath \cite{sagemath} for the lengthy symbolic computations. The resulting functions are explicit expressions in terms of $G$ and its derivatives, but for readability we only write out one of them here as an example:
\begin{align*}
    (\Delta_{p_1}F)(x,1,1,1,1,1,1,1)&=\frac{\phi_0x^4G'(x)^3}{G(x)^3} - \frac{\phi_0x^3G'(x)^2}{G(x)^2}.
\end{align*}
Now we can use the analytic information on the generating function $G(x)$ that Lemma~\ref{lem:bgw-gf} provides: analyticity in the $\Delta$-domain~\eqref{eq:delta-domain}
carries over from $G$ to its derivatives. Moreover, the equation $G(x) = x \Phi(G(x))$ ensures that $G(x)$ does not have any zeros other than $x = 0$ in its domain of analyticity. Thus we know now that $(\Delta_{p_1}F)(x,1,1,1,1,1,1,1)$ is analytic in the same $\Delta$-domain.

Now we recall from Lemma~\ref{lem:bgw-gf} that, as $x \to 1$,
\begin{equation}\label{G-singular-exp}
G(x)=1-\frac{\sqrt{2}}{\sigma}\sqrt{1-x}+o(|x-1|^{1/2}).    
\end{equation}
Such an asymptotic expansion carries over to derivatives of $G$, see \cite[Section VI.10.1]{Flajolet2009} for details (the proof is given for a big-$O$ error term, but it works equally well with little-$o$). Hence we have
\begin{equation}\label{G1-singular-exp}
G'(x)=\frac{\sqrt{2}}{2\sigma}(1-x)^{-1/2}+o(|1-x|^{-1/2}).
\end{equation}
This can be continued for higher derivatives.
So taking derivatives and plugging into the expressions for the different generating functions $(\Delta_{p_i}F)(x,1,1,1,1,1,1,1)$, we get
\begin{align*}
        (\Delta_{p_1}F)(x,1,1,1,1,1,1,1)&=\frac{\sqrt{2} \phi_{0}}{4 \, \sigma^{3}}  {(1-x)}^{-\frac{3}{2}} + o\!\left( {|x-1|^{-\frac{3}{2}}}\right),\\                                                                                   
    (\Delta_{p_2}F)(x,1,1,1,1,1,1,1)&=\frac{5\sqrt{2}}{8\sigma^3}\left(-\frac{ \phi_0^{2}}{2} + \frac{ \phi_0}{2} + \frac{\phi_0^{2}}{\, \sigma^{2}}\right)  {(1-x)}^{-\frac{7}{2}} + o\left( {|x-1|}^{-\frac{7}{2}}\right),\\
     (\Delta_{p_3}F)(x,1,1,1,1,1,1,1)&=\frac{15\sqrt{2}\phi_0^2}{32\sigma^3}(1-x)^{-\frac{7}{2}}+o\left(|1-x|^{-\frac{7}{2}}\right),\\
      (\Delta_{p_4}F)(x,1,1,1,1,1,1,1)&=\frac{9009\sqrt{2}}{512\sigma^3}\left(-\frac{\phi_0^4}{2}+\frac{\phi_0^3}{2}+\frac{7\phi_0^4}{\sigma^2}\right)(1-x)^{-\frac{15}{2}}+o\left(|1-x|^{-\frac{15}{2}}\right).
\end{align*}
We can now use the Flajolet--Odlyzko singularity analysis \cite{Flajolet1990} (see also the books \cite{Flajolet2009} and \cite{Melczer2021}).
This yields asymptotic formulas for the $n$-th coefficients of the generating functions, which we then divide by the $n$-th coefficient of $G$ (whose asymptotics are obtained in the same way) to get the statement of the lemma.
\end{proof}

\begin{example}\label{example_GW}
We illustrate the general procedure by determining an asymptotic formula for the expectation $\E(|T|\pt-\lt\p)$.
We start by calculating the first order derivative with respect to $y$. Differentiate both sides of \eqref{functionalExpression} with respect to $y$ and evaluate both sides at $y=z=u=v=w=1$ to get
\[F_y(x)=x\Phi'(G(x))\cdot F_y(x)+x\phi_0,\]
where $F_y(x):=\left(\frac{\partial}{\partial y}F(x,y,z,u,v,w,s,t)\right)\big|_{y=z=u=v=w=s=t=1}$. In this calculation we also used that
\[\left(\frac{\partial}{\partial y}F(xu,yz,z,u,v,w,w,v)\right)\big|_{y=z=u=v=w=1}= \left(\frac{\partial}{\partial y}F(x,y,z,u,v,w,s,t)\right)\big|_{y=z=u=v=w=s=t=1}.\]
Solving this for $F_y(x)$, we get
\[F_y(x)=\frac{x\phi_0}{1-x\Phi'(G(x))}.\]
To continue, we return to \eqref{g=xphig}, which we differentiate and rewrite to get
\[1-x\Phi'(G(x))=\frac{G(x)}{xG'(x)}.\]
This can then be plugged into the expression for $F_y(x)$ to get
\[F_y(x)=\frac{\phi_0 x^2 G'(x)}{G(x)}.\]
For the derivative in $z$, with $F_z(x)$ being defined analogously to $F_y(x)$, we obtain
\[F_z(x)=x\Phi'(G(x))\left(F_z(x)+F_y(x)\right),\]
since by the chain rule,
\[\left(\frac{\partial}{\partial z}F(xu,yz,z,u,v,w,w,v)\right)\big|_{y=z=u=v=w=1}=F_z(x)+F_y(x).\]
Solving for $F_z(x)$ and plugging in the expressions for $F_y(x)$ and $\Phi'(G(x))$ that have already been obtained, we get
\[F_z(x)=\frac{\phi_0x^3G'(x)^2}{G(x)^2}-\frac{\phi_0x^2G'(x)}{G(x)}.\]

Now we consider higher-order derivatives, starting with
\[F_{xy}(x):=\left(\frac{\partial^2}{\partial x \partial y}F(x,y,z,u,v,w,s,t)\right)\big|_{y=z=u=v=w=s=t=1}.\]
The generating functions $F_{xz}(x)$ and $F_{yu}(x)$ are defined in an analogous way.
Differentiating first with respect to $y$ and then with respect to $x$ using the chain rule produces an expression for $F_{xy}(x)$ involving a combination of derivatives of $\Phi$ and first order derivatives that have already been found.
\[F_{xy}(x)=\Phi'(G(x))F_y(x)+x\Phi''(G(x))F_x(x)F_y(x)+x\Phi'(G(x))F_{xy}(x) + \phi_0.\]
The new term $\Phi''(G(x))$ is calculated by differentiating \eqref{g=xphig} twice, which gives
\[\Phi''(G(x))=\frac{xG(x)G''(x)+2G(x)G'(x)-2xG'(x)^2}{x^3G'(x)^3}.\] We have $F_x(x)=G'(x)$ by definition, and an expression for $F_y(x)$ has already been determined. We can therefore again solve for $F_{xy}(x)$ and get
\[F_{xy}(x)=\frac{2 \phi_0 x G'(x)}{G(x)}+\frac{\phi_0 x^2 G''(x)}{G(x)}-\frac{\phi_0 x^2 G'(x)^2}{G(x)^2}.\]
Moving on to $F_{xz}(x)$, the chain rule gives
\[F_{xz}(x)=\Phi'(G(x))(F_z+F_y)+x\Phi''(G(x))F_x(x)(F_z(x)+F_y(x))+x\Phi'(G(x))(F_{xz}(x)+F_{xy}(x)).\]
Now all parts (with the exception of $F_{xz}(x)$) have been handled already, and so we can plug these previous results in and solve for $F_{xz}(x)$ toget
\[F_{xz}(x)=\phi_0\left(\frac{4x^2 G'(x)^2}{G(x)^2}-\frac{2x^3 G'(x)^3}{G(x)^3}+\frac{2x^3 G'(x)G''(x)}{G(x)^2}-\frac{2x G'(x)}{G(x)}-\frac{x^2 G''(x)}{G(x)}\right).\]
In a similar fashion, we obtain
\[F_{yu}(x)=\phi_0\left(\frac{5x^3G'(x)^2}{G(x)^2}-\frac{3x^4G'(x)^3}{G(x)^3}+\frac{2x^4G'(x)G''(x)}{G(x)^2}-\frac{2x^2G'(x)}{G(x)}-\frac{x^3G''(x)}{G(x)}\right).\]
We can now conclude that 
\begin{align*}
    (\Delta_{p_1}F)(x,1,1,1,1,1,1,1)&=xF_{xz}(x)-F_{yu}(x)=\frac{\phi_0x^4G'(x)^3}{G(x)^3} - \frac{\phi_0x^3G'(x)^2}{G(x)^2}.
\end{align*}
Plugging in~\eqref{G-singular-exp},~\eqref{G1-singular-exp} and the corresponding expansion of $G''(x)$ gives us
\begin{align*}
 (\Delta_{p_1}F)(x,1,1,1,1,1,1,1)&=\frac{\sqrt{2} \phi_{0}}{4 \, \sigma^{3}}  {(1-x)}^{-\frac{3}{2}} + o\!\left( {|x-1|^{-\frac{3}{2}}}\right).
\end{align*}
Now singularity analysis yields
\[[x^n](\Delta_{p_1}F)(x,1,1,1,1,1,1,1)\sim \frac{\sqrt{2}\phi_0}{4\sigma^3}\frac{n^{1/2}}{\Gamma(3/2)}=\frac{\phi_0}{\sigma^3\sqrt{2\pi}}n^{1/2}\]
and
\[[x^n]G(x)\sim -\frac{\sqrt{2}}{\sigma}\frac{n^{-3/2}}{\Gamma(-1/2)}=
\frac{1}{\sigma\sqrt{2\pi}}n^{-3/2}.\]
So we can finally conclude that
\[\E(|T|\pt-\lt\p) = \frac{[x^n](\Delta_{p_1}F)(x,1,1,1,1,1,1,1)}{[x^n]G(x)} \sim \frac{\phi_0}{\sigma^2} n^2.\]
\end{example}

We now have all the necessary ingredients to prove Theorem~\ref{Theorem_simplygenerated}. We will apply Lemma~\ref{XoverYLemma} to the quotient of the random variables $|T|P_E(T) - |L(T)|P(T)$ and $|L(T)|$, as well as $\bt\wt-\blt\w$ and $\bt\blt$, and use the asymptotic formulas in Lemma~\ref{Lemma_simplygenerated}. The required concentration property of the number of leaves is provided by Lemma~\ref{lem:concentration-bgw}.

\begin{proof}[Proof of Theorem~\ref{Theorem_simplygenerated}]
We use Lemma~\ref{lem:concentration-bgw} with $x = n^{3/4}$, which gives us
\[ \P(|\lt-\phi_0n| > n^{3/4}) \leq \kappa \sqrt{n} e^{-\lambda \sqrt{n}}.\]
Thus, we can apply Lemma~\ref{XoverYLemma}, combined with the first part of Lemma~\ref{Lemma_simplygenerated}, to
$X_n=|T|\pt-\lt\p$ and $Y_n=|T|\lt$ with $f(n) = \kappa \sqrt{n} e^{-\lambda \sqrt{n}}$, $g(n) = n^{7/4}$, $h(n) = n^3$ and $\mu_n = \phi_0n^2$. Boundedness of $\frac{\E(X_n)}{\mu_n}$ from below is clear, and boundedness of $\frac{\E(X_n)}{\E(|X_n|)}$ from below follows from the Cauchy--Schwarz inequality and Lemma~\ref{Lemma_simplygenerated}.
This readily 
proves the first statement \eqref{path1Thm} of the theorem. 
For the remaining parts of the theorem, we first use Lemma~\ref{Lemma_highermomentsconsentration} with the bounds $\lt\leq n$ and $\phi_0n\leq n$ to get
\begin{equation}\label{eq:leaf_square_concentration}
    \P(|\lt^2-\phi_0^2n^2| > 2n^{7/4})\leq \kappa \sqrt{n} e^{-\lambda \sqrt{n}}.
\end{equation}
This allows us to apply Lemma~\ref{XoverYLemma} with appropriately chosen $X_n$ and $Y_n$ together with the second part of Lemma~\ref{Lemma_simplygenerated} to finish the proof of \eqref{TheoremEqSecondMomentGW} in Theorem~\ref{Theorem_simplygenerated}. Equations \eqref{path1Thm} and \eqref{TheoremEqSecondMomentGW} show that the variance of the difference tends to a constant, and Chebyshev's inequality yields boundedness in probability.
\end{proof}
The proof of Theorem~\ref{Theorem:BGW_Wiener} follows the same structure.
\begin{proof}[Proof of Theorem~\ref{Theorem:BGW_Wiener}]
    Since $\blt=\frac{1}{2}(\lt^2-\lt)$ and $\lt\leq n$, we have
    \begin{equation}\label{eq:binomial_concentration}
        \left|\blt-\frac{\phi_0^2n^2}{2}\right|\leq \frac{1}{2}\left|\lt^2-\phi_0^2n^2\right|+\frac{\lt}{2}\leq n^{7/4}+\frac{n}{2}\leq 2n^{7/4}
    \end{equation}
    whenever $|\lt^2-\phi_0^2n^2| \leq 2n^{7/4}$. Therefore, by~\eqref{eq:leaf_square_concentration},
    \begin{equation*}
        \P\left( \left|\blt-\frac{\phi_0^2n^2}{2}\right|> 2 n^{7/4}\right)\leq \kappa \sqrt{n} e^{-\lambda \sqrt{n}}.
    \end{equation*}
    Applying Lemma~\ref{XoverYLemma} together with the third part of Lemma~\ref{Lemma_simplygenerated} to the random variables $X_n=\bt\wt-\blt\w$ and $Y_n=\bt\blt$ with $f(n) = \kappa \sqrt{n} e^{-\lambda \sqrt{n}}$, $g(n) = n^{15/4}$, $h(n) = n^5$ and $\mu_n = \binom{n}{2}\frac{\phi_0^2 n^2}{2}$ proves the first statement of the theorem. 
    
    The only part that does not follow as in the proof of Theorem~\ref{Theorem_simplygenerated} is the required lower bound of $Y_n$, which however is easily resolved as the only problematic case where there is only one leaf occurs with exponentially small probability (compare Remark~\ref{remark:BGW_tree_is_path}). We can condition on $\lt>1$ to force $Y_n \geq 1$, or simply define $Y_n = 1$ if $\lt = 1$. The impact is exponentially small in any case. 

    For the second part of the theorem, Lemma~\ref{Lemma_highermomentsconsentration} together with $\blt\leq \frac{n^2}{2}$ and $\frac{\phi_0^2n^2}{2}\leq \frac{n^2}{2}$ gives
    \begin{equation*}
         \P\left( \left|\blt^2-\frac{\phi_0^4n^4}{4}\right|>  2 n^{15/4}\right)\leq \kappa \sqrt{n} e^{-\lambda \sqrt{n}}.
    \end{equation*}
    Applying Lemma~\ref{XoverYLemma} with the last part of Lemma~\ref{Lemma_simplygenerated} together with suitable choices of $X_n$ and $Y_n$ finishes the proof of~\eqref{eq:BGW_Wiener_Thm_eq2} in Theorem~\ref{Theorem:BGW_Wiener}. Together, \eqref{eq:BGW_Wiener_Thm_eq1} and~\eqref{eq:BGW_Wiener_Thm_eq2} show that the variance of the difference tends to a constant, and by Chebyshev's inequality the difference is again bounded in probability.
\end{proof}
\begin{proof}[Proof of Corollary~\ref{cor:airy}]
  Since
  \[\frac{\pt}{\lt} - \frac{\p}{|T|}\]
  is bounded in probability,
  \[n^{-1/2} \Big( \frac{\pt}{\lt} - \frac{\p}{|T|} \Big) = \frac{\pt}{\sqrt{n} \lt} - \frac{\p}{n^{3/2}}\]
  goes to $0$ in probability. It follows from Slutsky's theorem that $\frac{\pt}{\sqrt{n} \lt}$ converges weakly to the same distribution as $\frac{\p}{n^{3/2}}$, which we know to be $\sigma^{-1} \eta$. Since $\frac{\lt}{n}$ converges in probability to $\phi_0$ (a consequence of Lemma~\ref{lem:concentration-bgw}), a second application of Slutsky's theorem proves the statement.
\end{proof}

\begin{proof}[Proof of Corollary~\ref{cor:BGW_Wiener}]
    The proof proceeds exactly as the proof of Corollary~\ref{cor:airy}. The convergence of $\frac{\w}{n^{5/2}}$ in distribution to $\sigma^{-1}\zeta$ was proved by Janson in \cite{Janson03}.
\end{proof}

\section{Increasing trees}\label{sec:LPA}

We proceed along the same lines as in the previous section.
The main lemma needed to prove Theorem~\ref{sats_alltogether} and Theorem~\ref{Theorem:increasing_Wiener} is

\begin{lemma}\label{Lemma_alltogether}
Let $T$ be as in Theorem~\ref{sats_alltogether}. Then we have
\begin{equation*}
    \E\left(|T|\pt-\lt\p\right)\sim \frac{1}{(\alpha+2)^2}n^2
\end{equation*}
and
\begin{equation*}
   \E\left(\left(|T|\pt-\lt\p\right)^2\right)\sim\frac{1}{(\alpha+2)^4}n^4
\end{equation*}
as $n \to \infty$. Similarly,
\begin{equation*}
\E\left(\binom{|T|}{2}\wt-\binom{\lt}{2}\w\right)\sim \frac{1}{2}\frac{(\alpha+1)}{(\alpha+2)^3}n^4 
\end{equation*}
and
\begin{equation*}
    \E\left(\left(\binom{|T|}{2}\wt-\binom{\lt}{2} \w\right)^2\right)\sim \frac{1}{4}\frac{(\alpha+1)^2}{(\alpha+2)^6} n^8 
\end{equation*}
as $n\to\infty$.
\end{lemma}

\begin{proof}
Similarly to the Bienaym\'e--Galton--Watson case, we aim to extend (\ref{recursiveODE}), (\ref{daryODE}) and (\ref{prefattODE}) to include all of the variables of interest. In general, we construct the exponential generating function 
\[F(x,y,z,u,v,w,s,t)=\sum_T \omega(T) \frac{1}{|T|!}x^{|T|}y^{\lt} z^{\pt} u^{\p}v^{\qt}w^{\q}s^{|T|^2}t^{\lt^2}, \]
where the sum is over all rooted ordered increasing trees. For all three families, we can use the same recursions as before to get a differential equation
\begin{equation*}
\frac{\partial}{\partial x}F(x,y,z,u,v,w,1,1)=\varphi{\left(F(xu,yz,z,u,v,w,w,v)\right)}+y-1,
\end{equation*}
where the correction term $y-1$ accounts for the case of the single vertex tree and $\varphi(\cdot)$ is a function depending on the family. In particular, we get
\begin{align}
    \frac{\partial}{\partial x}\Fr(x,y,z,u,v,w,1,1)&=\exp \left(\Fr(xu,yz,z,u,v,w,w,v)\right)+y-1,\label{recursivePDE}\\
    \frac{\partial}{\partial x}\Fd(x,y,z,u,v,w,1,1)&={\left(1+\Fd(xu,yz,z,u,v,w,w,v)\right)^d}+y-1,\nonumber\\
    \frac{\partial}{\partial x}\Fp(x,y,z,u,v,w,1,1)&={\left(1-\Fp(xu,yz,z,u,v,w,w,v)\right)^{-\frac{1}{\alpha}}}+y-1.\nonumber
\end{align}
The subscript of $F$ indicates the type of tree, with $\text{r}$ indicating recursive trees, $d$ indicating $d$-ary trees and $\alpha$ indicating GPORTs.
Using the same notation as previously, we get that
\begin{align*}
    (\Delta_{p_1}F)(x,1,1,1,1,1,1,1)&=\sum_T \omega(T)\big(|T|\pt-\lt\p\big) \frac{x^{|T|}}{|T|!},\\
    (\Delta_{p_2}F)(x,1,1,1,1,1,1,1)&=\sum_T \omega(T)\big(|T|\pt-\lt\p\big)^2 \frac{x^{|T|}}{|T|!}
\end{align*}
and
\begin{align*}
    (\Delta_{p_3}F)(x,1,1,1,1,1,1,1)&=\sum_T \omega(T)\left(\binom{|T|}{2}\wt-\binom{\lt}{2}\w\right)\frac{x^{|T|}}{|T|!},\\
    (\Delta_{p_4}F)(x,1,1,1,1,1,1,1)&=\sum_T \omega(T)\left(\binom{|T|}{2}\wt-\binom{\lt}{2}\w\right)^2\frac{x^{|T|}}{|T|!}.
\end{align*}

Dividing the $n$-th coefficient in the respective generating functions by the corresponding $n$-th coefficient in $G$ yields the expected values of interest. Note here that the factor $\frac{1}{|T|!}$ is automatically cancelled. To illustrate the calculations, Example~\ref{example_recursive_leaves} considers the expectation of $|T|\pt-\lt\p$ in a random recursive tree.
Finding other derivatives of $F$ is done in a similar manner, which always involves solving an ODE of the form
\[f'(x)=a(x)f(x)+b(x),\]
where $a(x)$ and $b(x)$ are functions dependent on lower order derivatives of $F$ resulting from the iterative calculation of inner derivatives and applications of the product rule. In fact, we always have $a(x)=\varphi'(G(x))$. Such an ODE has the general solution
\[f(x)=Ce^{A(x)}+e^{A(x)}\int b(x)e^{-A(x)} dx,\]
where $A(x)$ is any antiderivative of $a(x)$.
Carrying out these calculations\footnote{Once again, all code is available at the following repository:
\url{https://github.com/JacobLundblad/ExternalPathLengthAndTerminalWienerIndex}.} in SageMath \cite{sagemath}, we get explicit functions in all these cases, but the expressions are quite large and also require some extra case distinctions. 
Using error terms to increase readability, we have
\begin{align*}
    (\Delta_{p_1}\Fr)(x,1,1,1,1,1,1,1)&=\frac{1}{4}(1-x)^{-2}+O(\log(1-x)),\\
    (\Delta_{p_2}\Fr)(x,1,1,1,1,1,1,1)&=\frac{3}{8}(1-x)^{-4}+O(\log(1-x)(1-x)^{-3}),\\
      (\Delta_{p_3}\Fr)(x,1,1,1,1,1,1,1)&=\frac{3}{8}(1-x)^{-4}+O(\log(1-x)(1-x)^{-3}),\\
        (\Delta_{p_4}\Fr)(x,1,1,1,1,1,1,1)&=\frac{315}{16}(1-x)^{-8}+O(\log(1-x)(1-x)^{-7}),
\end{align*}
for recursive trees. Similarly,
\begin{align*}
    (\Delta_{p_1}\Fd)(x,1,1,1,1,1,1,1)&=\frac{(\beta - 1)^2(\beta)_2}{(\beta - 2)^2} t^{\beta-2}+O(t^{\beta-1}\log t ),\\
    (\Delta_{p_2}\Fd)(x,1,1,1,1,1,1,1)&=\frac{(\beta - 1)^4(\beta)_4}{(\beta - 2)^4} t^{\beta-4}+O(t^{\beta-3}\log t ),\\
      (\Delta_{p_3}\Fd)(x,1,1,1,1,1,1,1)&=-\frac{1}{2}\frac{(\beta-1)^2 (\beta)_4}{(\beta-2)^3}t^{\beta-4}+O(t^{\beta-3}\log t),\\
        (\Delta_{p_4}\Fd)(x,1,1,1,1,1,1,1)&=\frac{1}{4}\frac{(\beta-1)^4(\beta)_8}{(\beta-2)^6}t^{\beta-8}+O(t^{\beta-7}\log t),
\end{align*}
for $d$-ary trees with $\beta=\frac{1}{1-d}$ and $t=1-(d-1)x$, where $(\cdot)_k$ denotes the falling factorial.
In the final case of GPORTs, we get, with $\beta=\frac{\alpha}{1+\alpha}$ and $t=1-(1+\frac{1}{\alpha})x$,
\begin{align*}
    (\Delta_{p_1}\Fp)(x,1,1,1,1,1,1,1)&=-\frac{(\beta - 1)^2(\beta)_2}{(\beta - 2)^2} t^{\beta-2}+O(t^{\beta-1}\log t ),\\
    (\Delta_{p_2}\Fp)(x,1,1,1,1,1,1,1)&=-\frac{(\beta - 1)^4(\beta)_4}{(\beta - 2)^4} t^{\beta-4}+O(t^{\beta-3}\log t ),\\
     (\Delta_{p_3}\Fp)(x,1,1,1,1,1,1,1)&=\frac{1}{2}\frac{(\beta-1)^2 (\beta)_4}{(\beta-2)^3}t^{\beta-4}+O(t^{\beta-3}\log t),\\
      (\Delta_{p_4}\Fp)(x,1,1,1,1,1,1,1)&=-\frac{1}{4}\frac{(\beta-1)^4(\beta)_8}{(\beta-2)^6}t^{\beta-8}+O(t^{\beta-7}\log t).
\end{align*}
Consider the last case as an example. These generating functions have coefficients
\begin{align*}
   [x^n] (\Delta_{p_1}\Fp)(x,1,1,1,1,1,1,1)&\sim-\frac{(\beta - 1)^2(\beta)_2}{(\beta - 2)^2}(1+1/\alpha)^n\frac{n^{-\beta+1}}{\Gamma(-\beta+2)},\\
   [x^n] (\Delta_{p_2}\Fp)(x,1,1,1,1,1,1,1)&\sim-\frac{(\beta - 1)^4(\beta)_4}{(\beta - 2)^4}(1+1/\alpha)^n\frac{n^{-\beta+3}}{\Gamma(-\beta+4)},\\
     [x^n] (\Delta_{p_3}\Fp)(x,1,1,1,1,1,1,1)&\sim \frac{1}{2}\frac{(\beta-1)^2 (\beta)_4}{(\beta-2)^3}(1+1/\alpha)^n\frac{n^{-\beta+3}}{\Gamma(-\beta+4)},\\
       [x^n] (\Delta_{p_4}\Fp)(x,1,1,1,1,1,1,1)&\sim -\frac{1}{4}\frac{(\beta-1)^4(\beta)_8}{(\beta-2)^6}(1+1/\alpha)^n \frac{n^{-\beta+7}}{\Gamma(-\beta+8)}.
\end{align*}
Dividing by
\[[x^n]\Gp(x)\sim -(1+1/\alpha)^n\frac{n^{-\beta-1}}{\Gamma(-\beta)}\]
proves the lemma, using the fact that 
$\Gamma(-\beta+k)=(-1)^k (\beta)_k\Gamma(-\beta)$ for any positive integer $k$. At the end, everything can be expressed in terms of $\alpha$ as $\beta = \frac{\alpha}{1+\alpha}$. Similar calculations also apply to recursive trees (with $\alpha = \beta = 0$) and $d$-ary increasing trees (with $\alpha = - \frac{1}{d}$ and $\beta = \frac{1}{1-d}$).
\end{proof}

\begin{example}\label{example_recursive_leaves}
 Again, we illustrate the general procedure with an example. Suppose we aim to calculate the expectation $\E(|T|\pt-\lt\p)$ for a random recursive tree. To do this, we first differentiate both sides of (\ref{recursivePDE}) with respect to $y$ and evaluate both sides at $y=z=u=v=w=1$ to get
\[\frac{\partial}{\partial x}F_y(x)=1+F_y(x)\cdot\exp(\Gr(x)),\]
where 
\begin{align*}
F_y(x) &:=\left(\frac{\partial}{\partial y}\Fr(x,y,z,u,v,w,s,t)\right)\big|_{y=z=u=v=w=s=t=1} \\
&= \left(\frac{\partial}{\partial y}\Fr(xu,yz,z,u,v,w,w,v)\right)\big|_{y=z=u=v=w=1}.
\end{align*}
Plugging in $\Gr(x)=\log\left(\frac{1}{1-x}\right)$, we get the first order linear ODE 
\[F_y'(x)=1+F_y(x)\cdot (1-x)^{-1},\]
which is solved by
\[F_y(x)=C(1-x)^{-1}+(1-x)^{-1}\left(x-\frac{x^2}{2}\right).\]
To find the constant $C$, we consider first that $F_y(0)=C$. Since the constant coefficient of $\Fr$, and thus of $F_y(x)$, is $0$, 
we have $C=0$ and thus
\begin{align*}
    F_y(x)&=(1-x)^{-1}\left(x-\frac{x^2}{2}\right)\\
    &=\frac{1}{2}(1-x)^{-1}-\frac{1}{2}(1-x).
\end{align*}
Next, we calculate the analogously defined $F_z(x)$. With the same chain rule corrections as in Example~\ref{example_GW} we get (again using $\exp(\Gr(x))=(1-x)^{-1}$) the ODE
\[F'_z(x)=(F_y(x)+F_z(x))(1-x)^{-1}.\]
The initial condition of the ODE is again $F_z(0)=0$, and plugging in the solution for $F_y(x)$ we get
\[F_z(x)=-\frac{1}{4}(1-x)^{-1}-\frac{1}{2}(1-x)^{-1}\log(1-x)+\frac{1}{4}(1-x).\]
From this we obtain $F_{xz}(x)$ by differentiating $F_z(x)$, so
\[F_{xz}(x)=\frac{1}{4(1-x)^2}+\frac{-\log(1-x)}{2(1-x)^2}-\frac{1}{4}.\]
Similarly,
\[F_u'(x)=(xF_x(x)+F_u(x))(1-x)^{-1}\]
gives
\[F_u(x)=\frac{-\log(1-x)-x}{1-x}.\]
This is then used in
\[F'_{yu}(x)=\left((xF_x(x)+F_u(x))F_y(x)+xF_{xy}(x)+F_{yu}(x)\right)(1-x)^{-1}\]
together with $F_{xy}(x)=F'_y(x)$ to get
\[F_{yu}(x)=\frac{-\log(1-x)}{2(1-x)^2}+\frac{\log(1-x)}{2(1-x)}-\frac{\log(1-x)}{2}-\frac{1}{4(1-x)}+\frac{1-x}{4}.\]
Thus,
\[(\Delta_{p_1}\Fr)(x,1,1,1,1,1,1,1)=xF_{xz}(x)-F_{yu}(x)=\frac{1}{4(1-x)^2}-\frac{1}{4}+\frac{1}{2}\log(1-x).\]
The coefficients of this generating function are exactly 
\[\frac{n+1}{4}-\frac{1}{2n}\] for $n\geq 1$, but in general the asymptotic behaviour of the coefficients is determined by the dominating term around the singularity of the function, in this case the term $\frac{1}{4}(1-x)^{-2}$. Using singularity analysis, we get that
\[[x^n]\Fr(x,1,1,1,1,1,1,1)\sim \frac{1}{4}\frac{n^{1}}{\Gamma(2)}=\frac{1}{4}n.\]
To get the expectation, we then have to multiply by $n$ (since $[x^n] \Gr(x) = \frac1{n}$), which yields the result
\[\E(|T|\pt-\lt\p)\sim\frac{n^2}{4}.\]   
\end{example}
\color{black}

Now the proof of Theorem~\ref{sats_alltogether} follows exactly the same lines as the proofs of Theorem~\ref{Theorem_simplygenerated}, with Lemma~\ref{Lemma_alltogether} and Lemma~\ref{lem:concentration-lpat} playing the roles of Lemma~\ref{Lemma_simplygenerated} and Lemma~\ref{lem:concentration-bgw}, respectively.
Likewise, the proof of Theorem~\ref{Theorem:increasing_Wiener} follows the proof of Theorem~\ref{Theorem:BGW_Wiener}, using Lemma~\ref{Lemma_alltogether} and Lemma~\ref{lem:concentration-lpat} instead of Lemma~\ref{Lemma_simplygenerated} and Lemma~\ref{lem:concentration-bgw}.

\begin{proof}[Proof of Corollary~\ref{cor:Xalpha}]
  We can argue in a similar fashion as in the proof of Corollary~\ref{cor:airy}: by Lemma~\ref{Lemma_alltogether},
\[\frac{n \pt - \lt \p}{n^2} \xrightarrow{P} \frac{1}{(\alpha+2)^2}\]
and
\[\E \Big( \frac{n \pt - \lt \p}{n^2} \Big) \to \frac{1}{(\alpha+2)^2}\]
as $n \to \infty$. Consequently, the difference of the two, which can be written as
\[\frac{\pt - \E(\pt)}{n} - \frac{\lt \p - \E(\lt \p)}{n^2},\]
converges in probability to $0$. So by Slutsky's theorem, it suffices to show that the random variable $\frac{\lt \p - \E(\lt \p)}{n^2}$ has the desired limiting distribution. We split it into three parts:
\begin{align}
\frac{\lt \p - \E(\lt \p)}{n^2} &= \frac{\lt (\p - \E (\p))}{n^2}\nonumber \\
&\quad + \frac{(\lt - \E(\lt)) \E(\p)}{n^2} \label{eq:three-terms} \\
&\quad+ \frac{\E(\lt) \E(\p) - \E(\lt \p)}{n^2}. \nonumber
\end{align}
Now, we know that $\frac{\p - \E(\p)}{n}$ has limiting distribution $X_{\alpha}$, and $\frac{\lt}{n}$ goes to $\frac{1+\alpha}{2+\alpha}$ in probability by the concentration property in 
Lemma~\ref{lem:concentration-lpat}. So it only remains to show that the second and third term go to $0$ in probability, so that another application of Slutsky's theorem yields the statement:
\begin{itemize}
\item Observe that $\frac{\lt - \E(\lt)}{n^{1/2+\epsilon}} \xrightarrow{P} 0$ for every $\epsilon > 0$ by Lemma~\ref{lem:concentration-lpat}. Since we know that $\E(\p) = O(n \log n)$, the second term on the right-hand side of~\eqref{eq:three-terms} goes to $0$ in probability.
\item Finally, one can show by a direct generating function computation that $\E(\lt) \E(\p) - \E(\lt \p) = o(n^2)$. Alternatively, one can apply Lemma~\ref{lem:concentration-lpat}: since the event that $|\lt - \E(\lt)| > n^{3/4}$ has a probability that goes to zero faster than any power of $n$ while $\lt$ and $\p$ are polynomially bounded, one can condition on the event that $|\lt - \E(\lt)| \leq n^{3/4}$ at the expense of a negligible error term to find that
\begin{align*}
\E(\lt) \E(\p) - \E(\lt \p) &= \E \Big( \big(\E(\lt) - \lt\big) \p \Big) \\
&= O\big(n^{3/4} \E(\p)\big) = O(n^{7/4}\log n).
\end{align*}
This shows that the last term in~\eqref{eq:three-terms} also goes to $0$.
\end{itemize}
Putting everything together, we obtain the statement of the corollary.
\end{proof}

\begin{proof}[Proof of Corollary~\ref{cor:increasing_Wiener}]
    Following the same argument as in the proof of Corollary~\ref{cor:Xalpha}, the random variable $\frac{\wt-\E(\wt)}{n^2}$ has the same limiting distribution as 
    \begin{align*}
        &\frac{2\blt\left(\w-\E\left(\w\right)\right)}{n^4}\\
        &+\frac{2}{n^4}\left(\blt-\E\left(\blt\right)\right)\E\left(\w\right)\\
        &+\frac{2}{n^4}\left(\E\left(\blt\right)\E\left(\w\right)-\E\left(\blt\w\right)\right).
    \end{align*}
    By the assumption that $\frac{\w-\E\left(\w\right)}{n^2}$ has limiting distribution $Y_\alpha$, and because $\frac{\lt}{n}$ converges to $\frac{1+\alpha}{2+\alpha}$ in probability by Lemma~\ref{lem:concentration-lpat}, two more applications of Slutsky's Theorem yield the statement of the corollary, provided that the second and third terms go to 0 in probability. This is achieved as follows:
    \begin{itemize}
        \item By the same calculation as in~\eqref{eq:binomial_concentration}, Lemma~\ref{lem:concentration-lpat} and Lemma~\ref{Lemma_highermomentsconsentration} yield $\frac{\blt-\E\left(\blt\right)}{n^{7/4+\epsilon}}\xrightarrow{P} 0$ for every $\epsilon>0$. We know that $\E\left(\w\right)=O(n^2\log n)$, and so the second term goes to 0 in probability.

        \item As in the proof of Corollary~\ref{cor:Xalpha}, a direct computation or another application of Lemma~\ref{lem:concentration-lpat} show that
        \[\E\left(\blt\right)\E\left(\w\right)-\E\left(\blt\w\right)=o(n^4),\]
        so the last term goes to $0$ in probability as well.
    \end{itemize}
\end{proof}

\section{Further extensions}

The methods used in this paper could be adapted to other tree functionals, provided that they have a similar recursive decomposition that allows for the extension of the univariate generating function and the corresponding functional or differential equation.

One may also expect similar phenomena for other families of random trees beyond those considered here. One such example is the family of unlabelled trees, which also has a generating function satisfying a functional equation and would therefore be amenable to the same type of multivariate analysis. Asymptotics of the moments of the Wiener index were obtained with this method in \cite{WAGNER20121502}.

\bibliography{refs}

\end{document}